\documentclass{amsart}
\theoremstyle{plain}
\newtheorem{Thm}{Theorem}

\newtheorem{Coro}[Thm]{Corollary}

\newtheorem{Lem}[Thm]{Lemma}

\begin{document}
\large
\title[Examples for Cross Curvature Flow]
{Examples for Cross Curvature Flow on 3-Manifolds}

\author{Li MA and Dezhong Chen}

\address{Department of mathematical sciences \\
Tsinghua university \\
Beijing 100084 \\
China}

\email{lma@math.tsinghua.edu.cn} \dedicatory{}
\date{May 13th, 2004}

\begin{abstract}
Recently, B.Chow and R.S.Hamilton \cite{CH} introduced the cross
curvature flow on 3-manifolds. In this paper, we analyze two
interesting examples for this new flow. One is on a square torus
bundle over a circle, and the other is on a $S^{2}$ bundle over a
circle. We show that the global flow exist in both cases. But on
the former the flow diverges at time infinity, and on the latter
the flow converges at time infinity.
\end{abstract}

 \maketitle

\begin{section}{Introduction}
  Recently, B.Chow and R.S.Hamilton \cite{CH} introduced the cross
curvature flow on 3-manifolds. This flow has a very nice property
that it has at least two invariant sets of Riemannian metrics. In
particular, the set of metrics of negative curvature is preserved.
So one may try to use this flow to study the Geometrization
Conjecture.

The aim of this paper is give two examples which show that the
cross curvature flow is very degenerate in some cases. We show
that the global flow may diverge at time infinity in a square
torus bundle.

 We now begin with
the definition of the cross curvature flow.

  Let $(M,g)$ be a 3-dimensional Riemannian manifold. The Einstein
  tensor is
  $$
  P_{ij}=R_{ij}-\frac{1}{2}Rg_{ij},
  $$
  where $R_{ij}$ is the Ricci tensor and $R$ is the scalar
  curvature of the metric $g$. Let
  \begin{equation}
  P^{ij}=g^{ik}g^{jl}R_{kl}-\frac{1}{2}Rg^{ij}.
  \end{equation}
  Then B.Chow and R.Hamilton \cite{CH} define the cross curvature tensor by
  $$
  h_{ij}=(\frac{\det P^{kl}}{\det g^{kl}})\cdot V_{ij},
  $$
  where $V_{ij}$ is the inverse of $P^{ij}$. Clearly from this definition
  and the matrix theory we know that
  the cross curvature can also be defined even if $P^{ij}$ has no
  inverse.
  The cross
  curvature flow is defined by
   \begin{equation}
  \frac{\partial}{\partial t}g_{ij}=2h_{ij},
   \end{equation}
  if the sectional curvature is negative or
   \begin{equation}
  \frac{\partial}{\partial t}g_{ij}=-2h_{ij}
  \end{equation}
  if the sectional curvature is positive. In some cases, B.Andrews
  \cite{A} proved the convergence results.

  We remark that the validity of
  definition of $h_{ij}$ for $P^{ij}$ being degenerate can also be seen from
  the following fact.
  Let $\mu_{ijk}$ be the volume form of $g$ with the normalization
  that $\mu_{123}=1$ and raise
  indices by $$\mu^{ijk}=g^{ip}g^{jq}g^{kr}\mu_{pqr}.$$
  Then the
  cross curvature can be written as
  $$
h_{ij}=\frac{1}{8}\mu^{pqk}\mu^{rsl}R_{ilpq}R_{kjrs}.
  $$
where $R_{ijkl}$ is the Riemannian curvature tensor of $g$.

  Roughly speaking, we show the following results:

\noindent(a) The cross
  curvature flow (XCF)(2)
has global solutions on a square torus bundle over a circle and
the global flow diverges at time infinity.

\noindent(b) The cross
  curvature flow (XCF)(3)
has global solutions on a $S^{2}$ bundle over a circle and the
global flow converges at time infinity.

For more precise assumptions on the initial metrics, one can
consult Theorems 1 and 11,14 in sections 3 and 4. We point out
that we reduce the XCF into degenerate parabolic equations. This
is the difficulty which does not appear in the Ricci-Hamilton flow
\cite{H95}. In section two, we recall some basic formulae in
Riemannian Geometry.

\end{section}

\begin{section}{Fundamental Local Formulas}
  It is convenient to compute using moving frame method. So, first of all, let us
  review some useful formulae.

  Let $M$ be a Riemannian manifold of dimension $n$.
  $\omega_{i},i=1,\cdots,n,$ is an orthogonal co-frame on $M$. $e_{i}$
  is its dual frame. Then we have the following structure equations:
  $$
  d\omega_{i}=\sum_{j}\omega_{j}\wedge\omega_{ji},
  $$
  $$
  d\omega_{ij}=\sum_{k}\omega_{ik}\wedge\omega_{kj}+\Omega_{ij}.
  $$
  The forms $\omega_{ij}$, which are anti-symmetric in $i,j$, are
  the connection forms. The $\Omega_{ij}$, also anti-symmetric in
  $i,j$, are the curvature forms. They have the expression
  \begin{equation}
  \Omega_{ij}=-\frac{1}{2}\sum_{k,l}R_{ijkl}\omega_{k}\wedge\omega_{l},
  \end{equation}
  where the coefficients $R_{ijkl}$ satisfy the well-known
  symmetry relations
  $$
  R_{ijkl}=-R_{jikl}=-R_{ijlk},
  $$
  $$
  R_{ijkl}=R_{klij},
  $$
  $$
  R_{ijkl}+R_{iklj}+R_{iljk}=0.
  $$

  The Ricci tensor $R_{ij}$ is defined by
  \begin{equation}
  R_{ij}=\sum_{k}R_{ikjk},
  \end{equation}
  and the scalar curvature by
  \begin{equation}
  R=\sum_{i}R_{ii}.
  \end{equation}

  In follows, we take $n=3$.
\end{section}

\begin{section}{XCF on a square torus bundle over a circle}
  Consider a 3-manifold $M^{3}$ where the torus group $T^{2}=S^{1}\times
  S^{1}$ acts freely. Then $M^{3}$ is a $T^{2}$ bundle over the
  circle $S^{1}$. There is a larger group $G$ which is the
  isometry group of the square flat torus $R^{2}/Z^{2}$,
  containing $T^{2}$ as a subgroup. Consider metrics on $M^{3}$
  which have $G$ as their isometry group with the subgroup $T^{2}$
  acting freely. These 3-manifolds are called square torus bundles
  over a circle. For details, see \cite{H95}. We now consider
  the metric on a square torus bundle taking the form
  $$
  ds^{2}=f(x)^{2}dx^{2}+g(x)^{2}(dy^{2}+dz^{2}),
  $$
  where $x$ is a coordinate on the orbit space $S^{1}$ and, $y$
  and $z$ are coordinates on the fibre, $f, g$ are positive
  smooth functions on $S^1$. Note that $ds=f(x)dx$ is
  the arc length for the quotient metric on the orbit space
  $S^{1}$, and $g(x)$ is the length of the side of the square
  fibre over $x$. For this metric, we can take the orthogonal
  co-frame $\omega_{i}$ as
  $$
  \omega_{1}=fdx,
  \omega_{2}=gdy,
  \omega_{3}=gdz.
  $$
  Then its dual frame is
  $$
  e_{1}=\frac{1}{f}\frac{\partial}{\partial x},
  e_{2}=\frac{1}{g}\frac{\partial}{\partial y},
  e_{3}=\frac{1}{g}\frac{\partial}{\partial z}.
  $$

  A straightforward computation shows that the connection forms
  are
  $$
  \omega_{12}=\frac{g_{x}}{f}dy,
  \omega_{13}=\frac{g_{x}}{f}dz,
  \omega_{23}=0,
  $$
  and the curvature forms are
  $$
  \Omega_{12}=(\frac{g_{x}}{f})_{x}dx\wedge dy,
  \Omega_{13}=(\frac{g_{x}}{f})_{x}dx\wedge dz,
  \Omega_{23}=(\frac{g_{x}}{f})^{2}dy\wedge dz.
  $$

  Using relation (4), we get
  $$
  R_{1212}=R_{1313}=-\frac{1}{fg}(\frac{g_{x}}{f})_{x},
  R_{2323}=-\frac{1}{g^{2}}(\frac{g_{x}}{f})^{2}.
  $$
 Other components of $R_{ijkl}$ are zeros.

  Substituting these expressions into equation (5), we get the Ricci
  tensor:
  $$
  R_{11}=-\frac{2}{fg}(\frac{g_{x}}{f})_{x},
  R_{22}=R_{33}=-\frac{1}{fg}(\frac{g_{x}}{f})_{x}-\frac{1}{g^{2}}(\frac{g_{x}}{f})^{2}.
  $$

  Substituting these expression into equation (6), we get the scalar
  curvature:
  $$
  R=-\frac{4}{fg}(\frac{g_{x}}{f})_{x}-\frac{2}{g^{2}}(\frac{g_{x}}{f})^{2}.
  $$

  Now by equation (1), we have
  $$
  P^{11}=\frac{1}{g^{2}}(\frac{g_{x}}{f})^{2},
  P^{22}=P^{33}=\frac{1}{fg}(\frac{g_{x}}{f})_{x}.
  $$
  Then we can compute the cross curvature tensor:
  $$
  h_{11}=[\frac{1}{fg}(\frac{g_{x}}{f})_{x}]^{2},
  h_{22}=h_{33}=\frac{1}{fg^{3}}(\frac{g_{x}}{f})^{2}(\frac{g_{x}}{f})_{x}.
  $$

  For later use, we should express the cross curvature tensor in
  natural coordinates $(x,y,z)$:
  $$
  h(\frac{\partial}{\partial x},\frac{\partial}{\partial x})
  =f^{2}h_{11}=[\frac{1}{g}(\frac{g_{x}}{f})_{x}]^{2},
  $$
  $$
  h(\frac{\partial}{\partial y},\frac{\partial}{\partial y})
  =g^{2}h_{22}=\frac{1}{fg}(\frac{g_{x}}{f})^{2}(\frac{g_{x}}{f})_{x},
  $$
  $$
  h(\frac{\partial}{\partial z},\frac{\partial}{\partial z})
  =g^{2}h_{33}=\frac{1}{fg}(\frac{g_{x}}{f})^{2}(\frac{g_{x}}{f})_{x}.
  $$

  So the XCF (2)
  $$
  \frac{\partial}{\partial t}g_{ij}=2h_{ij}
  $$
  reduces to the systems of degenerate evolution equations
$$
\{
    \begin{array}{ll}
    \frac{\partial f}{\partial t}=\frac{1}{fg^{2}}[(\frac{g_{x}}{f})_{x}]^{2},\\
    \frac{\partial g}{\partial t}=\frac{1}{fg^{2}}(\frac{g_{x}}{f})^{2}(\frac{g_{x}}{f})_{x}.\\
    \end{array}
$$

  Introduce the unit vector field on the orbit space
  $$
  \frac{\partial}{\partial s}=\frac{1}{f}\frac{\partial}{\partial x}
  $$
  whose evolution is given by the commutator
  $$
  [\frac{\partial}{\partial t},\frac{\partial}{\partial s}]
  =-\frac{1}{f}\frac{\partial f}{\partial t}\frac{\partial}{\partial s}.
  $$
  Then the XCF takes the form of the fully nonlinear degenerate parabolic equation
  $$
  \frac{\partial g}{\partial t}
  =\frac{1}{g^{2}}(\frac{\partial g}{\partial s})^{2}\frac{\partial^{2}g}{\partial s^{2}}
  $$
  on a circle whose unit vector field $\frac{\partial}{\partial s}$
  varies by the commutator
  $$
  [\frac{\partial}{\partial t},\frac{\partial}{\partial s}]
  =-\frac{1}{g^{2}}(\frac{\partial^{2}g}{\partial s^{2}})^{2}\frac{\partial}{\partial s}.
  $$
The local existence of the flow can be obtained from an
approximate method. In fact, for $\epsilon>0$, we consider
$$
  \frac{\partial g}{\partial t}
  =\frac{1}{g^{2}}[(\frac{\partial g}{\partial s})^{2}+\epsilon]\frac{\partial^{2}g}{\partial s^{2}}
  $$
We know from the standard parabolic equation theory that there is
a local solution for every $\epsilon>0$. Using the argument below
(see Lemmata 2,3,and 4) we can show that we have a global solution
for every $\epsilon>0$ with uniform bounds of derivatives.  Hence,
by sending $\epsilon\to 0$, we can prove that the cross curvature
flow (2) has a global solution for any initial data. Hence we have
\begin{Thm}
  The XCF on a square torus bundle over a circle has a
  solution which exists for all time.
  \end{Thm}

  Let
  $$
  g_{max}(t)=\max_{\theta\in S^{1}}g(\theta,t),
  $$
  and
  $$
  g_{min}(t)=\min_{\theta\in S^{1}}g(\theta,t).
  $$
  Throughout this section, we will assume that
  $$
  g_{max}(0)>g_{min}(0).
  $$
  \begin{Lem}
  For all $t>0$, $g_{max}(t)=g_{max}(0)$ and
  $g_{min}(t)=g_{min}(0)$.
  \end{Lem}
  \begin{proof}
  At the local maximum of $g$,
  $$
  \frac{\partial g}{\partial s}=0,
  $$
  so
  $$
  \frac{\partial g}{\partial t}=0,
  $$
 for $t>0$. The similar argument can be done for the local minimum of $g$.
 \end{proof}
We give here a remark. Using the commutator relation we have
  $$
  \frac{\partial}{\partial t}\frac{\partial g}{\partial s}
  =\frac{1}{g^{2}}(\frac{\partial g}{\partial s})^{2}\frac{\partial^{2}}{\partial s^{2}}\frac{\partial g}{\partial s}
  +\frac{1}{g^{2}}\frac{\partial g}{\partial s}(\frac{\partial}{\partial s}\frac{\partial g}{\partial s})^{2}
  -\frac{2}{g^{3}}(\frac{\partial g}{\partial s})^{3}\frac{\partial}{\partial s}\frac{\partial g}{\partial
  s}.
  $$
 Assume that the zero points of $g_s$ at $t=0$
 consist of  maximum or minimum points of $g$
 on $S^1$. We know from the equations above and the
 maximum principle (see Proposition 1.2 in \cite{An}) that the number of
 points with $g_s=0$
 is invariant along the flow.

  \begin{Lem}
  There exists a constant $C<\infty$ such that for all $t>0$,
  $$
  |\frac{\partial g}{\partial s}|\leq C.
  $$
  \end{Lem}
  \begin{proof}
 Recall that
  $$
  \frac{\partial}{\partial t}\frac{\partial g}{\partial s}
  =\frac{1}{g^{2}}(\frac{\partial g}{\partial s})^{2}\frac{\partial^{2}}{\partial s^{2}}\frac{\partial g}{\partial s}
  +\frac{1}{g^{2}}\frac{\partial g}{\partial s}(\frac{\partial}{\partial s}\frac{\partial g}{\partial s})^{2}
  -\frac{2}{g^{3}}(\frac{\partial g}{\partial s})^{3}\frac{\partial}{\partial s}\frac{\partial g}{\partial
  s}.
  $$

  This shows by the maximum principle
  that the maximum of $\frac{\partial g}{\partial s}$
  decreases and the minimum increases.
  \end{proof}
  \begin{Lem}
  There exists a constant $C<\infty$ such that for all $t>0$,
  $$
  |\frac{\partial^{2}g}{\partial s^{2}}|\leq C\exp(Ct).
  $$
  \end{Lem}
  \begin{proof}
  Repeatedly using the commutator relation we get
  $$
  \frac{\partial}{\partial t}(\frac{\partial^{2}g}{\partial s^{2}})^{2}
  =\frac{1}{g^{2}}(\frac{\partial g}{\partial s})^{2}\frac{\partial^{2}}{\partial s^{2}}
  (\frac{\partial^{2}g}{\partial s^{2}})^{2}-\frac{2}{g^{2}}(\frac{\partial g}{\partial s})^{2}
  (\frac{\partial^{3}g}{\partial s^{3}})^{2}-\frac{4}{g^{3}}(\frac{\partial g}{\partial s})^{3}
  \frac{\partial}{\partial s}(\frac{\partial^{2}g}{\partial s^{2}})^{2}
  $$
  $$
  \hspace{16mm}+\frac{4}{g^{2}}\frac{\partial g}{\partial s}\frac{\partial^{2}g}{\partial s^{2}}
  \frac{\partial}{\partial s}(\frac{\partial^{2}g}{\partial s^{2}})^{2}+[\frac{12}{g^{4}}
  (\frac{\partial g}{\partial s})^{4}-16\frac{\partial}{\partial t}\log(g)]
  (\frac{\partial^{2}g}{\partial s^{2}})^{2}.
  $$

  At the maximum of $(\frac{\partial^{2}g}{\partial s^{2}})^{2}$,
  $$
  \frac{\partial}{\partial s}(\frac{\partial^{2}g}{\partial s^{2}})^{2}=0,
  $$
  i.e.,
  $$
  \frac{\partial^{2}g}{\partial s^{2}}\frac{\partial^{3}g}{\partial s^{3}}=0.
  $$

  Then,
  $$
  \frac{\partial^{3}g}{\partial s^{3}}=0.
  $$
  Otherwise, $g$ will be constant which contracts Lemma 2.

  By Lemmas 2 and 3, we have that at the maximum of
  $(\frac{\partial^{2}g}{\partial s^{2}})^{2}$,
  $$
  \frac{\partial}{\partial t}\log(g^{16}(\frac{\partial^{2}g}{\partial s^{2}})^{2})\leq C
  $$
  for some constant $C<\infty$, from which we obtain the desired estimate.
  \end{proof}

  Now we make some geometrical observations about the XCF.
  \begin{Lem}
  The length $L$ of the orbit circle always increases.
  \end{Lem}
  \begin{proof}
  The arc length $ds$ on the orbit circle varies by
  $$
  \frac{\partial}{\partial t}ds=\frac{1}{g^{2}}(\frac{\partial^{2}g}{\partial s^{2}})^{2}ds
  $$
  and the length
  $$
  L=\int1ds
  $$
  varies by
  $$
  \frac{dL}{dt}=\int\frac{1}{g^{2}}(\frac{\partial^{2}g}{\partial s^{2}})^{2}ds\geq0.
  $$
  \end{proof}
  \begin{Thm}
  The length $L$ tends to infinity as $t\rightarrow\infty$.
  \end{Thm}
  \begin{proof}
  If not, then the time derivative of the length $L$ must be approaching zero at an $\epsilon$-dense
  set of sufficiently large times. If we can gain some control of its time derivative, then we can
  conclude that it converges to zero.
  $$
  \hspace{-60mm}\frac{d}{dt}\int\frac{1}{g^{2}}(\frac{\partial^{2}g}{\partial s^{2}})^{2}ds
  $$
  $$
  =\int-\frac{38}{3}\frac{1}{g^{5}}(\frac{\partial g}{\partial s})^{2}
  (\frac{\partial^{2}g}{\partial s^{2}})^{3}-\frac{1}{3}\frac{1}{g^{4}}
  (\frac{\partial^{2}g}{\partial s^{2}})^{4}-\frac{2}{g^{4}}(\frac{\partial g}{\partial s})^{2}
  (\frac{\partial^{3}g}{\partial s^{3}})^{2}
  $$
  $$
  \hspace{-50mm}+\frac{12}{g^{6}}(\frac{\partial g}{\partial s})^{4}
  (\frac{\partial^{2}g}{\partial s^{2}})^{2}
  $$
  $$
  \hspace{10mm}\leq\int(\frac{1}{\epsilon}+\frac{12}{g^{4}}(\frac{\partial g}{\partial s})^{4})
  \frac{1}{g^{2}}(\frac{\partial^{2}g}{\partial s^{2}})^{2}
  +\int(\epsilon\frac{361}{4}\frac{1}{g^{4}}(\frac{\partial g}{\partial s})^{4}-\frac{1}{3})
  \frac{1}{g^{4}}(\frac{\partial^{2}g}{\partial s^{2}})^{4}.
  $$
  We can choose $\epsilon$ small enough to make the integrand of the second term is negative.
  Then we have
  $$
  \frac{d}{dt}\int\frac{1}{g^{2}}(\frac{\partial^{2}g}{\partial s^{2}})^{2}ds
  \leq C\int\frac{1}{g^{2}}(\frac{\partial^{2}g}{\partial s^{2}})^{2}ds
  $$
  for some constant $C<\infty$. This bounds the growth of the integral to exponential. Hence it
  converges to zero as $t\rightarrow\infty$. This leads to
  $$
  \lim_{t\rightarrow\infty}\int(\frac{\partial^{2}g}{\partial s^{2}})^{2}ds=0.
  $$
  Then an easy argument shows that $g$ converges to a constant. This contracts Lemma 2.
  \end{proof}
  \begin{Coro}
  The total volume $V$ of the bundle increases to infinity.
  \end{Coro}
  \begin{proof}
  Since
  $$
  V=\int g^{2}ds
  $$
  we compute
  $$
  \frac{dV}{dt}=\int\frac{2}{3}\frac{1}{g^{2}}(\frac{\partial g}{\partial s})^{4}
  +(\frac{\partial^{2}g}{\partial s^{2}})^{2}\geq0.
  $$
  Furthermore,
  $$
  V\geq g_{min}(0)^{2}L.
  $$
  Then Lemma 2 and Theorem 6 show that $V\rightarrow\infty$ as $t\rightarrow\infty$.
  \end{proof}
\end{section}
\begin{section}{XCF on a $S^{2}$ bundle over a circle}
  Consider a 3-manifold $M^{3}$ where $S^{2}$ acts freely. Then
  $M^{3}$ is a $S^{2}$ bundle over the circle $S^{1}$. Topologically $M^{3}$ is
  $S^{2}\times S^{1}$. In this section,
  we consider the metric on $S^{2}\times S^{1}$
  which takes the form
  $$
  ds^{2}=f(x)^{2}dx+g(x)^{2}(dy^{2}+\cos^{2}ydz^{2}).
  $$
  Geometrically, $ds=f(x)dx$ is the arc length for the quotient metric
  on the orbit space
  $S^{1}$, and $g(x)$ is the length of the radius of the $S^{2}$
  fibre over $x$.For this metric, we can take the orthogonal
  co-frame $\omega_{i}$ as
  $$
  \omega_{1}=fdx,
  \omega_{2}=gdy,
  \omega_{3}=g\cos ydz.
  $$
  Then its dual frame is
  $$
  e_{1}=\frac{1}{f}\frac{\partial}{\partial x},
  e_{2}=\frac{1}{g}\frac{\partial}{\partial y},
  e_{3}=\frac{1}{g\cos y}\frac{\partial}{\partial z}.
  $$

  The connection forms are
  $$
  \omega_{12}=\frac{g_{x}}{f}dy,
  \omega_{13}=\frac{g_{x}}{f}\cos ydz,
  \omega_{23}=-\sin ydz.
  $$

  The curvature forms are
  $$
  \Omega_{12}=(\frac{g_{x}}{f})_{x}dx\wedge dy,
  \Omega_{13}=(\frac{g_{x}}{f})_{x}\cos ydx\wedge dz,
  $$
  $$
  \Omega_{23}=[(\frac{g_{x}}{f})^{2}-1]\cos ydy\wedge dz.
  $$

  Using equation (4), we get
  $$
  R_{1212}=R_{1313}=-\frac{1}{fg}(\frac{g_{x}}{f})_{x},
  R_{2323}=-\frac{1}{g^{2}}[(\frac{g_{x}}{f})^{2}-1].
  $$
 Other components are zeros.

  Then the Ricci tensor is
  $$
  R_{11}=-\frac{2}{fg}(\frac{g_{x}}{f})_{x},
  R_{22}=R_{33}=-\frac{1}{fg}(\frac{g_{x}}{f})_{x}-\frac{1}{g^{2}}[(\frac{g_{x}}{f})^{2}-1].
  $$
  The scalar curvature is
  $$
  R=-\frac{4}{fg}(\frac{g_{x}}{f})_{x}-\frac{2}{g^{2}}[(\frac{g_{x}}{f})^{2}-1].
  $$
  Now we can compute
  $$
  P^{11}=\frac{1}{g^{2}}[(\frac{g_{x}}{f})^{2}-1],
  P^{22}=P^{33}=\frac{1}{fg}(\frac{g_{x}}{f})_{x}.
  $$
  This leads to
  $$
  h_{11}=[\frac{1}{fg}(\frac{g_{x}}{f})_{x}]^{2},
  h_{22}=h_{33}=\frac{1}{fg^{3}}[(\frac{g_{x}}{f})^{2}-1](\frac{g_{x}}{f})_{x}.
  $$
  Still, we need the expression of the cross curvature tensor in
  natural coordinates:
  $$
  h(\frac{\partial}{\partial x},\frac{\partial}{\partial x})
  =f^{2}h_{11}=[\frac{1}{g}(\frac{g_{x}}{f})_{x}]^{2},
  $$
  $$
  h(\frac{\partial}{\partial y},\frac{\partial}{\partial y})
  =g^{2}h_{22}=\frac{1}{fg}[(\frac{g_{x}}{f})^{2}-1](\frac{g_{x}}{f})_{x},
  $$
  $$
  h(\frac{\partial}{\partial z},\frac{\partial}{\partial z})
  =g^{2}\cos^{2}yh_{33}=\frac{\cos^{2}y}{fg}[(\frac{g_{x}}{f})^{2}-1](\frac{g_{x}}{f})_{x}
  $$

  So the XCF
  $$
  \frac{\partial}{\partial t}g_{ij}=-2h_{ij}
  $$
  reduces to the systems of evolution equations
$$
\{
    \begin{array}{ll}
    \frac{\partial f}{\partial t}=-\frac{1}{fg^{2}}[(\frac{g_{x}}{f})_{x}]^{2},\\
    \frac{\partial g}{\partial t}=-\frac{1}{fg^{2}}[(\frac{g_{x}}{f})^{2}-1]
    (\frac{g_{x}}{f})_{x}.\\
    \end{array}
$$

  Define the unit vector field on the orbit space
  $$
  \frac{\partial}{\partial s}=\frac{1}{f}\frac{\partial}{\partial x}
  $$
  whose evolution is given by the commutator
  $$
  [\frac{\partial}{\partial t},\frac{\partial}{\partial s}]
  =-\frac{1}{f}\frac{\partial f}{\partial t}\frac{\partial}{\partial s}.
  $$
  Then the XCF takes the form of the fully nonlinear equation
  $$
  \frac{\partial g}{\partial t}
  =\frac{1}{g^{2}}[1-(\frac{\partial g}{\partial s})^{2}]\frac{\partial^{2}g}{\partial s^{2}}
  $$
  on a circle whose unit vector field $\frac{\partial}{\partial s}$ varies by the commutator
  $$
  [\frac{\partial}{\partial t},\frac{\partial}{\partial s}]
  =\frac{1}{g^{2}}(\frac{\partial^{2}g}{\partial s^{2}})^{2}\frac{\partial}{\partial s}.
  $$
Throughout this section, we will assume that the following
  initial condition holds
  $$
  \max_{S^{1}}|\frac{\partial g}{\partial s}(\cdot,0)|\leq\frac{1}{4}.
  $$
  Then the local solution exists by the standard parabolic
  equation theory.

  As before, applying the maximum principle to the evolution of $g$, we directly have
  \begin{Lem}
  $g_{max}(t)$ decreases and $g_{min}(t)$ increases.
  \end{Lem}

  \begin{Lem}
  For all $t>0$, we have
  $$
  \max_{S^{1}}|\frac{\partial g}{\partial s}(\cdot,t)|\leq\frac{1}{4}.
  $$
  \end{Lem}
  \begin{proof}
  Using the commutator relation
  $$
  \frac{\partial}{\partial t}(\frac{\partial g}{\partial s})^{2}
  =\frac{1}{g^{2}}(1-(\frac{\partial g}{\partial s})^{2})\frac{\partial^{2}}{\partial s^{2}}
  (\frac{\partial g}{\partial s})^{2}-\frac{2}{g^{3}}\frac{\partial g}{\partial s}
  (1-(\frac{\partial g}{\partial s})^{2})\frac{\partial}{\partial s}
  (\frac{\partial g}{\partial s})^{2}
  $$
  $$
  \hspace{-50mm}-\frac{2}{g^{2}}(\frac{\partial^{2}g}{\partial s^{2}})^{2}.
  $$
  The initial condition shows that
  $$
  \frac{\partial}{\partial t}(\frac{\partial g}{\partial
  s})^{2}\leq0
  $$
  at the maximum of $(\frac{\partial g}{\partial
  s})^{2}$. So the maximum decreases.
  \end{proof}
  \begin{Lem}
  There exists a constant $C<\infty$ such that
  $$
  |\frac{\partial^{2}g}{\partial s^{2}}|\leq C\exp(Ct)
  $$
  for all $t\geq0$.
  \end{Lem}
  \begin{proof}
  Repeatedly using the commutator relation
  \begin{equation}
  \frac{\partial}{\partial t}\frac{\partial^{2}g}{\partial s^{2}}
  =\frac{1}{g^{2}}(1-(\frac{\partial g}{\partial s})^{2})\frac{\partial^{2}}{\partial s^{2}}
  \frac{\partial^{2}g}{\partial s^{2}}+(\frac{4}{g^{3}}(\frac{\partial g}{\partial s})^{3}
  -\frac{4}{g^{3}}\frac{\partial g}{\partial s}-\frac{4}{g^{2}}\frac{\partial g}{\partial s}
  \frac{\partial^{2}g}{\partial s^{2}})\frac{\partial^{3}g}{\partial s^{3}}
  \end{equation}
  $$
  \hspace{4mm}+\frac{6}{g^{4}}(\frac{\partial g}{\partial s})^{2}(1-(\frac{\partial g}{\partial s})^{2})
  \frac{\partial^{2}g}{\partial s^{2}}-\frac{8}{g^{3}}(\frac{1}{4}-
  (\frac{\partial g}{\partial s})^{2})(\frac{\partial^{2}g}{\partial s^{2}})^{2}.
  $$

  If $g$ is constant, then the flow is over. So without loss of
  generality, we may assume $g$ is not constant. Then at the maximum
  of $\frac{\partial^{2}g}{\partial s^{2}}$,
  $$
  \frac{\partial^{2}g}{\partial s^{2}}>0,\frac{\partial^{3}g}{\partial
  s^{3}}=0,\frac{\partial^{4}g}{\partial s^{4}}\leq0.
  $$
  Note that the fourth term on the r.h.s. of equation (7) is non-positive. So
  $$
  \frac{\partial}{\partial t}\frac{\partial^{2}g}{\partial s^{2}}\leq
  C\frac{\partial^{2}g}{\partial s^{2}}
  $$
  for some $C>0$. Therefore
  \begin{equation}
  \frac{\partial^{2}g}{\partial s^{2}}\leq C\exp(Ct)
  \end{equation}
  for some $C>0$.

  Likewise, at the minimum of $\frac{\partial^{2}g}{\partial
  s^{2}}$,
  $$
  \frac{\partial^{2}g}{\partial s^{2}}<0,\frac{\partial^{3}g}{\partial
  s^{3}}=0,\frac{\partial^{4}g}{\partial s^{4}}\geq0.
  $$
  Then
  $$
  \frac{\partial}{\partial t}\log(-\frac{\partial^{2}g}{\partial s^{2}})
  \leq C-\frac{2}{g^{3}}(1-(\frac{\partial g}{\partial s})^{2})
  \frac{\partial^{2}g}{\partial s^{2}}+\frac{6}{g^{3}}(\frac{\partial g}{\partial s})^{2}
  \frac{\partial^{2}g}{\partial s^{2}}
  $$
  $$
  \hspace{-18mm}\leq C-\frac{\partial}{\partial t}\log g^{2}
  $$
  $\Rightarrow$
  $$
  \frac{\partial}{\partial t}\log(-g^{2}\frac{\partial^{2}g}{\partial
  s^{2}})\leq C
  $$
  $\Rightarrow$
  \begin{equation}
  \frac{\partial^{2}g}{\partial s^{2}}\geq-C\exp(Ct)
  \end{equation}
  for some $C>0$.

  Combining (8) and (9), we obtain the desired estimate.
  \end{proof}
  As a direct consequence, we have
  \begin{Thm} Assume $$
  \max_{S^{1}}|\frac{\partial g}{\partial s}(\cdot,0)|\leq\frac{1}{4}.
  $$
  The XCF on the $S^{2}$ bundle over a circle has a solution which
  exists for all time.
  \end{Thm}
  As before, we will give some geometrical results.
  \begin{Lem}
  The length $L$ of the orbit circle always decreases.
  \end{Lem}
  \begin{proof}
  The arc length $ds$ on the orbit circle varies by
  $$
  \frac{\partial}{\partial t}ds=-\frac{1}{g^{2}}(\frac{\partial^{2}g}{\partial
  s^{2}})^{2}ds.
  $$
  Then the length varies by
  $$
  \frac{dL}{dt}=-\int\frac{1}{g^{2}}(\frac{\partial^{2}g}{\partial
  s^{2}})^{2}ds\leq0.
  $$
  \end{proof}
  \begin{Lem}
  The $L^{2}$ norm of $\frac{\partial^{2}g}{\partial
  s^{2}}$ converges to zero as $t\rightarrow\infty$.
  \end{Lem}
  \begin{proof}
  Since we have known that $L$ decreases, so as in the proof of
  Theorem 6, what we need to do is to control
  $$
  \frac{d}{dt}\int\frac{1}{g^{2}}(\frac{\partial^{2}g}{\partial
  s^{2}})^{2}ds.
  $$
  A straightforward computation shows that
  $$
  \hspace{-76mm}\frac{d}{dt}\int\frac{1}{g^{2}}(\frac{\partial^{2}g}{\partial
  s^{2}})^{2}ds
  $$
  $$
  =\int-\frac{2}{g^{4}}(1-(\frac{\partial g}{\partial s})^{2})(\frac{\partial^{3}g}{\partial
  s^{3}})^{2}+\frac{1}{3g^{4}}(\frac{\partial^{2}g}{\partial
  s^{2}})^{4}+\frac{1}{g^{5}}(\frac{38}{3}(\frac{\partial g}{\partial s})^{2}-6)(\frac{\partial^{2}g}{\partial
  s^{2}})^{3}
  $$
  $$
  \hspace{-48mm}+\frac{12}{g^{6}}(\frac{\partial g}{\partial
  s})^{2}(1-(\frac{\partial g}{\partial s})^{2})(\frac{\partial^{2}g}{\partial
  s^{2}})^{2}.
  $$
  Note that
  $$
  -\frac{2}{g^{4}}(1-(\frac{\partial g}{\partial s})^{2})\geq-C_{1}
  $$
  for some $C_{1}>0$.

  Repeatedly using the $H\ddot{o}lder$ inequality and the
  Peter-Paul inequality yields
  $$
  \frac{d}{dt}\int\frac{1}{g^{2}}(\frac{\partial^{2}g}{\partial
  s^{2}})^{2}ds\leq-C_{1}\int(\frac{\partial^{3}g}{\partial
  s^{3}})^{2}ds+C_{2}(\sup(\frac{\partial^{2}g}{\partial
  s^{2}})^{2}+1)\int(\frac{\partial^{2}g}{\partial
  s^{2}})^{2}ds
  $$
  for some $C_{2}>0$.

  But
  $$
  \sup(\frac{\partial^{2}g}{\partial
  s^{2}})^{2}\leq\int|\frac{\partial}{\partial s}(\frac{\partial^{2}g}{\partial
  s^{2}})^{2}|ds
  $$
  $$
  \hspace{22mm}=2\int|\frac{\partial^{2}g}{\partial
  s^{2}}|\cdot|\frac{\partial^{3}g}{\partial
  s^{3}}|ds
  $$
  $$
  \hspace{35mm}\leq\int(\frac{\partial^{2}g}{\partial
  s^{2}})^{2}ds+\int(\frac{\partial^{3}g}{\partial
  s^{3}})^{2}ds.
  $$
  Therefore
  $$
  \hspace{-64mm}\frac{d}{dt}\int\frac{1}{g^{2}}(\frac{\partial^{2}g}{\partial
  s^{2}})^{2}ds
  $$
  $$\leq(C_{1}+C_{2})\int(\frac{\partial^{2}g}{\partial
  s^{2}})^{2}ds+\sup(\frac{\partial^{2}g}{\partial
  s^{2}})^{2}\cdot(C_{2}\int(\frac{\partial^{2}g}{\partial
  s^{2}})^{2}ds-C_{1})
  $$
  $$
  \leq C\int\frac{1}{g^{2}}(\frac{\partial^{2}g}{\partial
  s^{2}})^{2}ds+\sup(\frac{\partial^{2}g}{\partial
  s^{2}})^{2}\cdot(C_{0}\int(\frac{\partial^{2}g}{\partial
  s^{2}})^{2}ds-C_{1})
  $$
  for some positive constant $C$ and $C_{0}$. Here we have used Lemma 8. This bounds the
  growth of the integral to exponential when it is sufficiently
  small. Hence it converges to zero as $t\rightarrow\infty$.
  \end{proof}
  \begin{Thm} Assume $$
  \max_{S^{1}}|\frac{\partial g}{\partial s}(\cdot,0)|\leq\frac{1}{4}.
  $$
  Along the XCF on a $S^{2}$ bundle over a circle, $g$ converges as
  $t\rightarrow\infty$ to a constant $\alpha>0$.
  \end{Thm}
  \begin{proof}
  Since
  $$
  \sup(\frac{\partial g}{\partial
  s})^{2}\leq\int|\frac{\partial}{\partial s}(\frac{\partial g}{\partial
  s})^{2}|ds
  $$
  $$
  \hspace{22mm}=2\int|\frac{\partial g}{\partial
  s}|\cdot|\frac{\partial^{2}g}{\partial
  s^{2}}|ds
  $$
  $$
  \hspace{30mm}\leq C\cdot L^{\frac{1}{2}}\cdot(\int|\frac{\partial^{2}g}{\partial
  s^{2}}|^{2}ds)^\frac{1}{2}
  $$
  $$
  \hspace{-5mm}\rightarrow0
  $$
  as $t\rightarrow\infty$, so
  $$
  g_{max}-g_{min}\leq\int|\frac{\partial g}{\partial
  s}|ds
  $$
  $$
  \hspace{24mm}\leq\sup|\frac{\partial g}{\partial
  s}|\cdot L
  $$
  $$
  \hspace{8mm}\rightarrow0
  $$
  as $t\rightarrow\infty$.
  \end{proof}
  \begin{Lem}
  The $L^{2}$ norm of $\frac{\partial^{3}g}{\partial
  s^{3}}$ converges to zero as $t\rightarrow\infty$.
  \end{Lem}
  \begin{proof}
  We compute
  $$
  \hspace{-64mm}\frac{d}{dt}\int(\frac{\partial^{3}g}{\partial
  s^{3}})^{2}ds
  $$
  $$
  =\int-\frac{2}{g^{2}}(1-(\frac{\partial g}{\partial
  s})^{2})(\frac{\partial^{4}g}{\partial
  s^{4}})^{2}+\frac{24}{g^{4}}(\frac{\partial g}{\partial
  s})^{2}(1-(\frac{\partial g}{\partial
  s})^{2})(\frac{\partial^{3}g}{\partial
  s^{3}})^{2}
  $$
  $$
  -\frac{44}{g^{3}}(\frac{3}{11}-(\frac{\partial g}{\partial
  s})^{2})\frac{\partial^{2}g}{\partial
  s^{2}}(\frac{\partial^{3}g}{\partial
  s^{3}})^{2}+\frac{1}{g^{2}}(\frac{\partial^{2}g}{\partial
  s^{2}})^{2}(\frac{\partial^{3}g}{\partial
  s^{3}})^{2}
  $$
  $$
  \hspace{4mm}+\frac{8}{g^{2}}\frac{\partial g}{\partial
  s}\frac{\partial^{2}g}{\partial
  s^{2}}\frac{\partial^{3}g}{\partial
  s^{3}}\frac{\partial^{4}g}{\partial
  s^{4}}-\frac{120}{g^{6}}(\frac{\partial g}{\partial
  s})^{4}(1-(\frac{\partial g}{\partial
  s})^{2})(\frac{\partial^{2}g}{\partial
  s^{2}})^{2}
  $$
  $$
  \hspace{14mm}+\frac{248}{g^{5}}(\frac{\partial g}{\partial
  s})^{2}(\frac{15}{31}-(\frac{\partial g}{\partial
  s})^{2})(\frac{\partial^{2}g}{\partial
  s^{2}})^{3}-\frac{96}{g^{4}}(\frac{1}{8}-(\frac{\partial g}{\partial
  s})^{2})(\frac{\partial^{2}g}{\partial
  s^{2}})^{4}
  $$
  $$
  \hspace{-48mm}+\frac{32}{g^{3}}\frac{\partial g}{\partial
  s}(\frac{\partial^{2}g}{\partial
  s^{2}})^{3}\frac{\partial^{3}g}{\partial
  s^{3}}.
  $$

  When $t$ is sufficiently large, we have
  $$
  \frac{d}{dt}\int(\frac{\partial^{3}g}{\partial
  s^{3}})^{2}ds\leq-C_{1}\int(\frac{\partial^{4}g}{\partial
  s^{4}})^{2}ds+\sup(\frac{\partial^{3}g}{\partial
  s^{3}})^{2}+\int(\frac{\partial^{2}g}{\partial
  s^{2}})^{2}ds
  $$
  for some $C_{1}>0$.

  But
  $$
  \sup(\frac{\partial^{3}g}{\partial
  s^{3}})^{2}\leq\int|\frac{\partial}{\partial s}(\frac{\partial^{3}g}{\partial
  s^{3}})^{2}|ds
  $$
  $$
  \hspace{22mm}=2\int|\frac{\partial^{3}g}{\partial
  s^{3}}|\cdot|\frac{\partial^{4}g}{\partial
  s^{4}}|ds
  $$
  $$
  \hspace{46mm}\leq\frac{C_{1}}{2}\int(\frac{\partial^{4}g}{\partial
  s^{4}})^{2}ds+\frac{2}{C_{1}}\int(\frac{\partial^{3}g}{\partial
  s^{3}})^{2}ds,
  $$
  then
  $$
  \frac{d}{dt}\int(\frac{\partial^{3}g}{\partial
  s^{3}})^{2}ds\leq-\frac{C_{1}}{2}\int(\frac{\partial^{4}g}{\partial
  s^{4}})^{2}ds+\frac{2}{C_{1}}\int(\frac{\partial^{3}g}{\partial
  s^{3}})^{2}ds+\int(\frac{\partial^{2}g}{\partial
  s^{2}})^{2}ds.
  $$

  Note that
  $$
  \int(\frac{\partial^{3}g}{\partial
  s^{3}})^{2}ds=\int-\frac{\partial^{2}g}{\partial
  s^{2}}\frac{\partial^{4}g}{\partial
  s^{4}}ds\leq(\int(\frac{\partial^{2}g}{\partial
  s^{2}})^{2}ds\cdot\int(\frac{\partial^{4}g}{\partial
  s^{4}})^{2}ds)^{\frac{1}{2}}.
  $$
  So if we assume that
  $$
  \int(\frac{\partial^{3}g}{\partial
  s^{3}})^{2}ds>C_{0}\int(\frac{\partial^{2}g}{\partial
  s^{2}})^{2}ds
  $$
  for some $C_{0}>0$ to be determined, then we get
  $$
  \int(\frac{\partial^{3}g}{\partial
  s^{3}})^{2}ds\leq C^{-1}_{0}\int(\frac{\partial^{4}g}{\partial
  s^{4}})^{2}ds.
  $$
  Hence
  $$
  \frac{d}{dt}\int(\frac{\partial^{3}g}{\partial
  s^{3}})^{2}ds\leq-(\frac{C_{0}C_{1}}{2}-\frac{2}{C_{1}}-\frac{1}{C_{0}})
  \int(\frac{\partial^{3}g}{\partial
  s^{3}})^{2}ds.
  $$

  If we take $C_{0}$ large enough, for instance,
  $$
  C_{0}>\frac{1+\frac{2}{C_{1}}+\sqrt{(1+\frac{2}{C_{1}})^{2}+2C_{1}}}{C_{1}},
  $$
  then
  $$
  \frac{d}{dt}\int(\frac{\partial^{3}g}{\partial
  s^{3}})^{2}ds\leq-\int(\frac{\partial^{3}g}{\partial
  s^{3}})^{2}ds.
  $$
  Therefore, either $\int(\frac{\partial^{3}g}{\partial
  s^{3}})^{2}ds$ decays exponentially, or it is comparable to
  $\int(\frac{\partial^{3}g}{\partial
  s^{3}})^{2}ds$. In either event, it decreases to zero.
  \end{proof}

  Since
  $$
  \sup(\frac{\partial^{2}g}{\partial
  s^{2}})^{2}\leq\int|\frac{\partial}{\partial s}(\frac{\partial^{2}g}{\partial
  s^{2}})^{2}|ds
  $$
  $$
  \hspace{22mm}=2\int|\frac{\partial^{2}g}{\partial
  s^{2}}|\cdot|\frac{\partial^{3}g}{\partial
  s^{3}}|ds
  $$
  $$
  \hspace{44mm}\leq2(\int(\frac{\partial^{2}g}{\partial
  s^{2}})^{2}ds)^{\frac{1}{2}}\cdot(\int(\frac{\partial^{3}g}{\partial
  s^{3}})^{2}ds)^{\frac{1}{2}}
  $$
  $$
  \hspace{-6mm}\rightarrow0
  $$
  as $t\rightarrow\infty$, then we can see that
  $$
  K(\frac{\partial}{\partial x},\frac{\partial}{\partial y})
  =K(\frac{\partial}{\partial x},\frac{\partial}{\partial z})
  =-\frac{1}{g}\frac{\partial^{2}g}{\partial
  s^{2}}\rightarrow0,
  $$
  and
  $$
  K(\frac{\partial}{\partial y},\frac{\partial}{\partial z})
  =\frac{1}{g^{2}}(1-(\frac{\partial g}{\partial s})^{2})
  \rightarrow\alpha.
  $$
  For the scalar curvature we have
  $$
  R(t) \rightarrow2\alpha.
  $$
  \end{section}

\end{document}